# A DYNAMIC ANALYSIS OF SEGMENTED LABOR MARKET

**Patrice Gaubert[1] et Marie Cottrell[2]**

gaubert@univ-paris1.fr, cottrell@univ-paris1.fr

[1] MATISSE, UMR 8595 CNRS
MSE, Université Paris 1, 106, Bd de l'Hôpital, 75013 PARIS

[2] MATISSE-SAMOS, UMR 8595 CNRS
PMF, Université Paris 1, 90, rue de Tolbiac, 7634 PARIS CEDEX 13

# A DYNAMIC ANALYSIS OF SEGMENTED LABOR MARKET

Patrice Gaubert[1] et Marie Cottrell[2]

gaubert@univ-paris1.fr, cottrell@univ-paris1.fr

[1] MATISSE, UMR 8595 CNRS
MSE, Université Paris 1, 106, Bd de l'Hôpital, 75013 PARIS

[2] MATISSE-SAMOS, UMR 8595 CNRS
PMF, Université Paris 1, 90, rue de Tolbiac, 7634 PARIS CEDEX 13

**Abstract :**

*Using the Panel Study of Income Dynamics data on the period 1982-1992, this paper investigates some mechanisms of the labor market in the United States. This market is analyzed as a stable structure constituted of segments which present contrasted characteristics under the usual distinction between primary and secondary sectors. Using a neural network algorithm applied on quantitative variables measured at the level of heads of household, a broad classification in four classes of situations is constructed. It shows a clear hierarchy going from situations of very precarious work or no work at all, to situations of stable jobs with higher wages than the average. A Markov chain, constructed with the trajectories between the different situations of these workers, shows a very stable structure of this segmented labor market.*

**Theme**: Wage inequality and mobility, Internal labor markets and labor relations.

JEL-Code: J31, J41, J49

## 1. OBJECT OF THIS STUDY

The key idea of approaches of segmented labor markets[1] is that jobs have to be distinguished from a qualitative point of view, even if different distinctions not overlapping are used like internal – external markets or primary – secondary segments. The first one addresses the rules governing the relations between the workers and their occupation (hiring, career, exit from employment). The second one is based on a qualitative comparison of existing jobs (level of earnings, stability, career).

These two distinctions use numerous common factors but they introduce two different theoretical explanations of the labor market as well as their empirical validations.

Various recent examples may be found of this diversity: In a theoretical way may be cited the contributions proposed by Bulow & Summers (1986) in terms of efficiency wages or by Albrecht & Vroman (1992) in a more institutional approach. On the empirical ground one may cited the observation of a significant bias between the levels of earnings and the characteristic of stability of jobs situated in large firms, as in Oi (1990). Besides the existence of efficiency wages, a suggested explanation, is the presence of higher fixed costs induced by a specialized on the job training, or the practice of higher wages induced by a more constrained work associated to the discipline of a team production.

The simplest way to take into account a segmentation in an empirical manner is to use a subjective determination as in Theodossiou (1995) using the results of a survey: are belonging to the primary segment all the workers answering that they have a career in their job. This subjective information is controlled using a more objective variable, having received a specific on the job training more than one year. These elements, career, special training, are only partial characteristics of a segmentation. They do not constitute sufficient conditions of its existence and the method of identification seems to be a weak one.

The main objective here is to produce a robust representation of the segmentation of the labor market using both qualitative and quantitative variables to express its different dimensions. These variables have to be significant of the induced quality of a job, in terms of earnings (level and time variation), of duration (hours per week and weeks per year), and of seniority (in the current job compared to the number of years of presence on the labor market since the beginning of the working life). The combination of these quantitative variables may produce a qualitative evaluation of different jobs. This have to be completed by some qualitative characteristics like the type of occupation and the industry. Of course some personal information like the level of education has to be taken in account.

---

[1] see Doeringer & Piore (1971).

The final objective is not to obtain a representation of the segmentation, but this representation is a step to identify the trajectories followed by workers between these segments over a long period.

After a presentation of the data and the method used to identify a broad classification (section 2), the initial results are characterized in order to justify the final classification (section 3), then the main trajectories are presented and used to produce a measure of a stationary state representing the 'virtual' structure of the segmented labor market (section 4). Finally the main factors influencing the transitions between states characterizing the main trajectories are identified (section 5).

## 2. THE DATA

The data are obtained from the PSID, using the heads of households present in their family during the period 1982-1992. This panel contains much information about their job and also about their personal characteristics, At this time no more filtering has been used to select the individuals[2], even considering the over-representation of low-income black people. Among other results the classification produced is supposed to group this particular population.

Among the initial 4 000 heads of households obtained, only 2 507 are kept due to missing values concerning the major variables used. At the same time, in order to produce information on the evolution of some important characteristics of their job, the two first years are used only to construct rates of growth or variables in difference. So the trajectories are identified on the period 1984-1992.

In order to produce the classification 15 variables have been selected, 12 are observed and 3 constructed with observed variables. According to the methodology used by the authors of the panel data, personal variables are measured for the year of the survey, while data concerning the participation to the labor market are evaluated for the preceding year: the dimension of the family in 1992 is measured at the time of the survey and the number of work hours evaluated for this year is actually the number of hours worked in 1991. So the effective period of observation of the labor market is 1983-1991.

---

[2] except for those whose accuracy indices concerning major variables like earnings or work hours are worse than 2 in the PSID codification.

One may assume that there is an important distinction, concerning the quality of a job, to be considered along with the durability expressed by seniority in a specific place, that is the duration over each year expressed by the number of weeks worked and also the number of hours per week. This is a decisive characteristic to separate 'normal' jobs from various types of contingent work whose influence is growing at the end of the studied period[3].

Two personal variables used are:
- the age of head (AGEH)
- the family size (SIZFAM)

Several variables describe the major characteristics of the main job:
- the seniority in the current job (SENH)
- the number of years of work since the age of 18 (ANCH)
- the annual number of work hours (HMJH)
- the number of weeks worked in a year (WMJH)
- the number of hours per week (HWMJH)
- the hourly wage obtained (RSALH), expressed in real terms in order to eliminate the effect of inflation when evaluating the existence of a career[4]
- the annual rate of growth of the hourly wage which is computed using the two last years (GRSALH) and two other variables expressing a change in the duration of the main job, variation of the number of work hours per week (VHWMJH) and variation of the number of weeks worked in a year (VWMJH) since 1982
- two variables expressing a situation of non-work during the year, the number of weeks out of the labor force (WOUTH) and the number of weeks unemployed (WUNEH).

At last, two variables indicating that the worker is practicing one extra job or more
- the number of extra jobs (NBXJH)
- the corresponding annual number of work hours (HEXJH).[5]

Several qualitative variables are used in order to obtain a better illustration of the obtained classification:
- race of head of the household (white, black, other)[6]

---

[3] it would be of a greater importance for a study conducted on the 90's. This work is possible using the Early Releases of PSID data for the period 94-97, at this time it is only in progress and it is too early to indicate some results.
[4] using the implicit deflator of GDP.
[5] in the following tables and in the Kohonen maps the variables are in an alphabetical order: a table ordered of the variables is found in the appendix.
[6] At this time the gender characteristic is not used due to the fact that the head of the household is a woman only when no man is present in the family; a further study is to be conducted using the data, available in the PSID, concerning the wife involved in the labor market.

        - level of education
            . less than 5 grades
            . between 5 and 12 grades
            . high school achieved with terminal grade
            . professional or academic training until BA
            . BA and more
        - type of occupation
            . technical and professional workers
            . managers and administrators
            . clerks
            . craftsmen
            . operatives
            . laborers
            . others
        - sector of activity
            . agricultural, mining and construction
            . manufacturing
            . transportation and communication
            . wholesale and retail trade
            . finance, insurance and real estate
            . business and personal services
            . professional services
            . public administration

The frequencies of these attributes for the whole sample are presented together with those computed for the main classes obtained (table 5).

## 3. A NEURAL NETWORK CLASSIFICATION

In order to obtain a valid classification for the whole period, the data used to construct the neural network is constituted of the observations made for the 15 variables for the three years, 1984 1988 1992, as if they were different individuals observed at the same time. So the classification is conducted on 7521 observations (2 507 x 3) and 15 variables.

Each variable being measured in its specific units, the very different sizes of the numbers for these variables could disturb the classification, so all the variables are standardized before the construction of a Kohonen grid. The dimensions of this grid are 8 x 8.

It has been establish, Kohonen (1995), Cottrell, Fort & Pagès (1998), Kaski & Oja (1999), that a Kohonen classification is an extension of k-means method with an imposed structure of neighborhood. One of the improvements resulting is that observations which are close or similar are placed in the same class or in neighbor classes. The neighborhood structure selected here is a two-dimensional grid. It produces 64 classes, each one being represented by its code vector as in the following figure 1. It can be seen that close code vector have similar profiles.

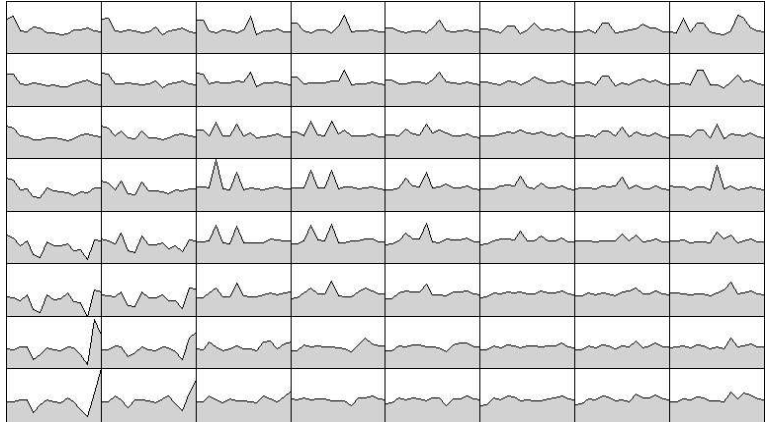

**Figure 1 : Profiles of the 64 code-vectors**

Each profile represent the vector defining a class, the values of each of the 15 variables being joined by a line for a better visualization. The classes are numbered from 1 to 64 when reading the grid from the upper left corner to the lower right one.

On the map obtained, where the variables are presented in each cell in alphabetical order, one may identify, in the lower left corner, classes containing individuals with no job (out of the labor force or unemployed) most of the year, in the central region, classes with people exerting more than one job at the same time, in the upper right corner, job situations with stability and high pay. The main diagonal represents a growing quality of the job situation, and the secondary one shows a clear opposition between the older workers in the upper left and the younger ones in the lower right.

With this classification giving 64 classes with the whole population observed over three years, it is possible to identify the complete trajectory followed by each individual from 1984 to 1992. This is obtained using the data available for each of the six other years, which are not used for the classification, to compute, for each individual the class of the closest code vector in terms of Euclidean distance[7].

For each individual the observed trajectory is represented by a sequence of 9 integers in the interval 1-64 and it is represented, according to the suggestion made by Serrano-Cinca in Deboeck & Kohonen (1998).

---

[7] when some values are missing this distance is computed with the existing values.

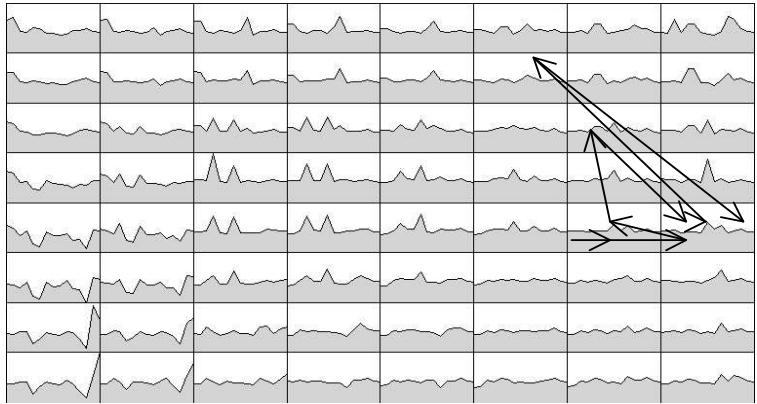

**Figure 2 : Trajectory of an individual staying in good job situation during the whole period**

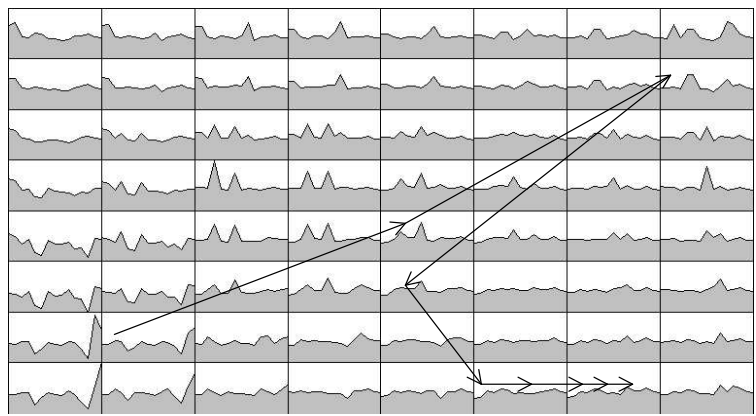

**Figure 3 : Trajectory of an individual leaving the more precarious situation to reach, after one year in a good situation, an intermediate position**

With 64 classes the maximum number of different trajectories is $64^9$ (eventually no more than the number of individuals in the sample). So it is very difficult to describe and interpret a great number of different trajectories, as it is not easy to interpret even the 64 classes obtained. This is why it is useful to reduce this number to something closer to the theoretical distinction of the labor market.

## 4. THE FINAL CLASSIFICATION (7 SUPER-CLASSES)

Different methods are used in order to reduce the number of classes (Cottrell & Rousset, 1997), for instance a hierarchical classification applied on the 64 code vectors. The result obtained is an unordered set of classes. Here the Kohonen algorithm has been used again, but the one dimensional version, producing an ordered sequence of broad classes which may be easily interpreted as going from the worst situation (unemployment and precarious jobs with low pay) to the better one (only one full time job with a great stability) via different intermediate situations (two or more extra jobs, partial time work, real wages under the average).

After several tries, the more significant structure seems to be a set of seven super-classes. The means of the 15 variables used to construct the Kohonen classification are given in table 1, for the whole sample and for each of the super-classes. For each variable the greatest value is in bold and the lowest in italics.

|        | Whole Population | Class 1 | Class 2 | Class 3 | Class 4 | Class 5 | Class 6 | Class 7 |
|--------|------------------|---------|---------|---------|---------|---------|---------|---------|
| AGEH   | 40.12 | 36.4 | *35.32* | **59.41** | *33.18* | 40.58 | **52.69** | 39.46 |
| ANCH   | 15.43 | *10.56* | 11.26 | **30.32** | *8.65* | 16.18 | **28.20** | 14.86 |
| GRSALH | 0.06 | *-0.18* | 0.02 | 0.07 | 0.06 | 0.03 | 0.02 | **0.19** |
| HEXJH  | 60.70 | 12.98 | **562.12** | 56.01 | *0.25* | 215.01 | *7.39* | *4.74* |
| HMJH   | 1974 | *663* | 1994 | *901* | 2040 | 2136 | 2008 | **2348** |
| HWMJH  | 42.18 | *24.69* | 41.88 | *22.95* | 42.09 | 44.34 | 42.09 | **48.72** |
| NBXJH  | 0.18 | *0.05* | **1.24** | 0.28 | *0* | 1.03 | *0.06* | *0.03* |
| RSALH  | 13.35 | *6.47* | 10.60 | 10.95 | 11.30 | 14.77 | 13.88 | **17.70** |
| SENH   | 91.14 | *19.51* | 64.02 | 41.05 | 58.28 | 118.81 | **173.39** | 93.04 |
| SIZFAM | 3.17 | 2.92 | 2.93 | *2.04* | 2.67 | **3.88** | 2.57 | **4.08** |
| VHWMJH | 0.59 | *-6.43* | 0.06 | *-17* | -0.13 | 0.23 | -0.52 | **5.23** |
| VWMJH  | 0.65 | *-15.66* | 2.77 | -3.83 | **3.89** | 0.17 | 1.05 | 2.92 |
| WMJH   | 44.61 | *15.29* | 47.51 | 40.81 | **48.48** | 48.23 | 47.60 | 48.10 |
| WOUTH  | 0.69 | **5.76** | *0.09* | 1.37 | *0.13* | *0.05* | *0.06* | *0.11* |
| WUNEH  | 2.09 | **16.08** | 0.80 | 3.29 | 0.40 | *0.13* | 0.41 | 0.53 |
| Size   | 7521 | 772 | 588 | 79 | 1932 | 416 | 1495 | 2240 |

**Table 1 : General mean and mean by super-class**

The different super-classes have very different sizes (See Figure 4), but the most important fact is that each one is clearly distinguished from the others using their quantitative characteristics.

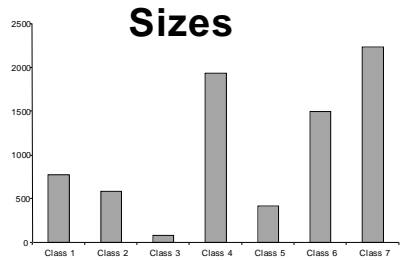

**Figure 4 : Sizes of the 7 super-classes.**

The limits of these super-classes have been represented on the Kohonen map (Figure 5), using different tones of gray. The Kohonen algorithm allows the visualization of the neighborhood between individuals and between classes.

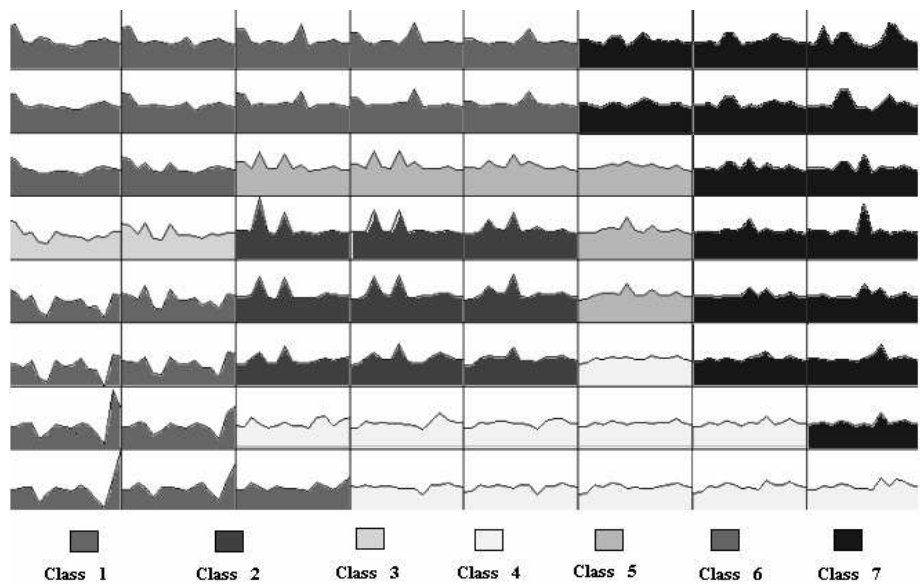

**Figure 5 : The 7 super-classes**

The main characteristics presented by the super-classes show the growing quality of the situation on the labor market, when moving from the lower left corner towards the upper right one. It is possible to summarize the super-classes as follows:

    - class 1 (lower left) is made of individuals younger than average population, with very short seniority in the current job, working a number of hours per week and a number of weeks per year neatly smaller than the average, with no extra job, low paid, being out of the labor force much more

often than the average and whose variations of hourly earnings and of hours of work over the period are mostly negative.

- class 2 (center of the map) contains individuals younger than the average, with a full-time job and one or more extra job, with earnings obtained in the main job severely lower than the average and no real changes in their situation over the period studied. This seems correspond to a situation where the main occupation, even a full-time one, produce an insufficient income and leads to the obligation to exert an extra job.

- class 3 (on top of class 1) is a very small one, approximately 50 persons each year, and presumably not really representative of a part of the whole population. Nevertheless this class is well defined with very old people having a kind of reduced activity which may be linked to a possible preparation to retirement.

- class 4 (below class 2) is made of young people with short seniority in their job, wages slightly below the average, an important augmentation in the number of hours worked but not in the wage obtained.

- class 5 (on top of class 2) is close from class 2 with respect to the fact that people have two or more occupations at the same time but very different when earnings are considered: it seems to represent in this class the possibility for skilled workers to obtain supplemental benefits from their special ability. This may be confirmed using information on the type of occupation and of industry.

- class 6 (upper left) appears to be made of older people with stable employment, having only one job and earnings that are close to the average. This class may be interpreted as a part of the primary segment, some lower tier or half, with respect to the presence at the same time of a great stability and a modest compensation.

- class 7 (upper right) is constituted of people of middle age, with large families, a good stability in their employment, working a longer duration than the average and for an hourly wage clearly above the average, with a steady growth of their compensation. It looks like the upper part of the primary segment.

This sequence of classes seems to be well ordered with respect to a growing quality of the situation occupied on the labor market. This may be seen with the graphical representation (figure 6) of the code vectors obtained for the super-classes with a one-dimensional algorithm.

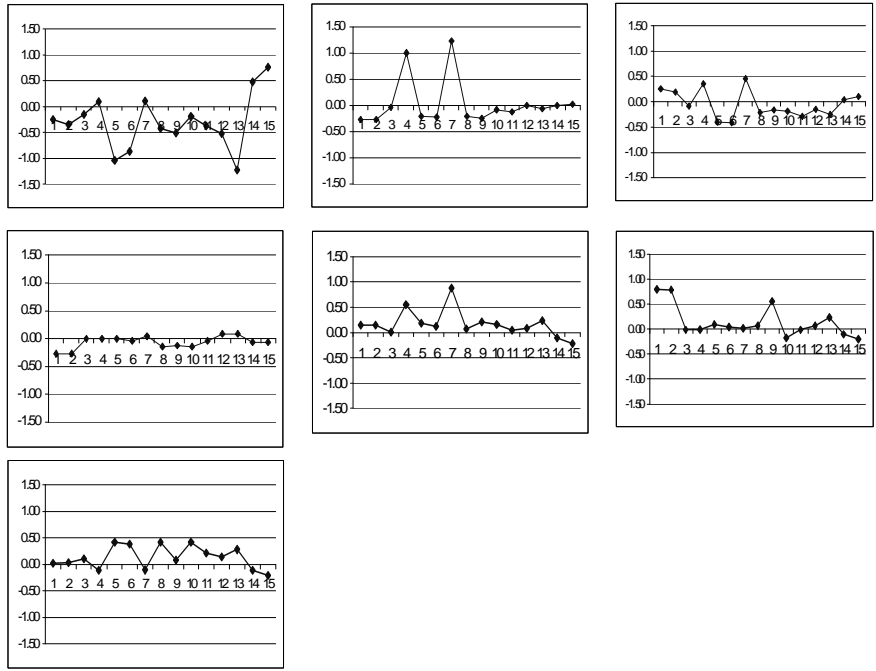

**Figure 6 : Code-vectors of the 7 super-classes**

The numbers appearing in the figure represent the 15 variables used ordered in alphabetical order (see appendix).

A first representation of the trajectories over the period of observation may be traced using these 7 classes, each individual being located each year in one of the class. A trajectory is represented by a sequence of 9 integers, as 234774344 for instance.

A transition is a sequence of two situations expressing that the individual has moved to another class. Considering the above example this individual is staying two consecutive years in the same situation and this two times, states 7 and 4, and is changing each year for the other years.
The diversity of situations may be summarized by presenting the frequencies of situations identifying the proportions of the different states occupied by an individual who is present most of the time in one specific class (5 years or more). They are in fact the probabilities to be in one other class when the dominant situation is a given class (table 2).

| Dominant Position | Size | Proba of 1 | Proba of 2 | Proba of 3 | Proba of 4 | Proba of 5 | Proba of 6 | Proba of 7 |
|---|---|---|---|---|---|---|---|---|
| 1 | 157 | **0.75** | 0.03 | 0.01 | 0.11 | 0.00 | 0.03 | 0.08 |
| 2 | 115 | 0.04 | **0.70** | 0.00 | 0.13 | 0.07 | 0.00 | 0.05 |
| 3 | 10 | 0.16 | 0.01 | **0.64** | 0.01 | 0.00 | 0.16 | 0.02 |
| 4 | 599 | 0.07 | 0.06 | 0.00 | **0.77** | 0.02 | 0.00 | 0.08 |
| 5 | 65 | 0.01 | 0.11 | 0.00 | 0.01 | **0.70** | 0.05 | 0.12 |
| 6 | 498 | 0.03 | 0.01 | 0.02 | 0.01 | 0.02 | **0.86** | 0.06 |
| 7 | 732 | 0.03 | 0.02 | 0.00 | 0.07 | 0.03 | 0.03 | **0.82** |

**Table 2 : Probabilities to be one year in a class being most of the time in a given class**

A group of people are not in a dominant position, that means that over 9 years no one of the 7 classes appears at least 5 times: 337 individuals (over 2507). For this group the probabilities are given in table 4.

| No dominant position | Size | Proba of 1 | Proba of 2 | Proba of 3 | Proba of 4 | Proba of 5 | Proba of 6 | Proba of 7 |
|---|---|---|---|---|---|---|---|---|
| | 331 | 0.14 | 0.16 | 0.03 | 0.22 | 0.13 | 0.08 | 0.23 |

**Table 3 : Probabilities when no class has a dominant position**

The main result is that the individuals stay most of the time in the same class, except for this group representing less than 15% of the whole sample. That means that the structure that appears constituted by segments with very different properties with respect to stability, existence of a career, seems to present the quality of a permanent state over a long period. This will be even clearer with the construction of a Markov chain.

Note that except for the very small group 3, the less stable classes, relatively, are those corresponding to lower situations, and precisely the two classes having extra job(s).

## 5. CLUSTERING INTO 4 MAIN CLASSES

At this stage, we could study the frequencies of the different trajectories, the frequencies of the category changes, the transition probabilities, etc; But the number of possibilities is still too large with 7 super-classes. So, to study more precisely the trajectories, a last grouping is realized, dividing the initial Kohonen map in four main classes :

    - super-classes 1 and 3 are grouped to produce the main class A: it is made of the more precarious conditions, recurring unemployment and low pay;

    - super-classes 2, 4 and 5 represent intermediate conditions: important duration of work and moderate wages. They constitute the main class B;

- super-classes 6 and 7 are still separated to produce the main classes C and D, respectively the lower part and the upper part of the primary segment.

A complimentary information is obtained with a principal component analysis concerning the role of the different variables used to construct the classification. The 5 first components 2/3 of the explanation. The projection of the variables on the first 3 axes are clearly significant (Appendix graphs 1-3).

The first axis is defined by the variables of activity: the number of work hours, the number of weeks, opposed to the number of weeks of unemployment and out of the labor force.

The second one oppose age, seniority to the family size (younger family are larger).

The third one is only defined with the extra job variables.

The level and the growth of wage and the variables in variation appears only as fourth and fifth axes. That means that the separation of the different situations is mainly explained by other factors than the differentiation of wages.

Even with this new grouping the main classes are well defined using the 15 quantitative variables (table 4), on the whole observations used to construct the initial classification (3 years pooled). The major characteristics observed above with the more detailed partition are still visible: work duration, seniority, level and growth of real wages, the practice of extra jobs. . For each variable the greatest value is in bold and the lowest in italics.

|        | Whole Population | Class A | Class B | Class C | Class D |
|--------|------------------|---------|---------|---------|---------|
| AGEH   | 40.12 | 38.56 | *34.66* | **52.69** | 39.46 |
| ANCH   | 15.43 | 12.39 | *10.24* | **28.20** | 14.86 |
| GRSALH | 0.06 | *-0.15* | 0.05 | 0.02 | **0.19** |
| HEXJH  | 60.70 | 16.97 | **143.25** | *7.39* | *4.74* |
| HMJH   | 1974 | *685.26* | 2045.05 | 2008.23 | **2348.84** |
| HWMJH  | 42.18 | *24.51* | 42.37 | 42.09 | **48.72** |
| NBXJH  | 0.18 | 0.07 | **0.40** | 0.06 | 0.03 |
| RSALH  | 13.35 | *6.88* | 11.66 | 13.88 | **17.70** |
| SENH   | 91.14 | *21.51* | 67.99 | **173.39** | 93.04 |
| SIZFAM | 3.17 | 2.85 | 2.89 | 2.57 | **4.08** |
| VHWMJH | 0.59 | *-7.41* | -0.04 | -0.52 | **5.23** |
| VWMJH  | 0.65 | *-14.57* | **3.14** | 1.06 | 2.92 |
| WMJH   | 44.61 | *17.66* | 48.25 | 47.60 | 48.10 |
| WOUTH  | 0.69 | **5.35** | 0.11 | *0.06* | 0.11 |
| WUNEH  | 2.09 | **14.89** | 0.44 | 0.41 | 0.53 |
| Size   | 7521 | 851 | 2936 | 1495 | 2240 |

**Table 4 : Means of the whole sample and by main class**

Using the Kohonen algorithm supply an order on the 4 situations A, B, C et D, from the more precarious and degraded to the more stable and with the best wages.

Some clear information may be obtained from 4 available qualitative variables, in terms of frequencies (see Table 5 for the data related to year 1992, and Appendix 3 for years 1988 and 1984).

| in 1992 | Whole sample | Class A | Class B | Class C | Class D |
|---|---|---|---|---|---|
| RACE | | | | | |
| 1 Whites | 69.1 % | 51.6 | 69.2 | 69.5 | 74.5 |
| 2 Blacks | 29.7 | 48.0 | 29.8 | 28.9 | 23.8 |
| EDUCATION | | | | | |
| 1 Primary | 0.9 | 1.4 | 0.2 | 2.6 | 0.4 |
| 2 Secondary | 18.9 | 32.7 | 13.1 | 27.8 | 15.1 |
| 3 Sec. achieved | 40.2 | 42.7 | 44.4 | 37.0 | 36.9 |
| 4 Post-sec. | 28.6 | 18.5 | 32.7 | 20.7 | 32.6 |
| 5 BA & more | 11.3 | 4.6 | 9.5 | 12.0 | 15.0 |
| OCCUPATION | | | | | |
| 0 No | 2.0 | 17.4 | 0 | 0 | 0.1 |
| 1-2 Managers, professionals | 36.2 | 15.7 | 36.7 | 32.9 | 44.9 |
| 4 Clerks | 12.0 | 14.2 | 13.2 | 14.0 | 8.8 |
| 5 Craftsmen | 17.1 | 14.9 | 17.0 | 18.5 | 17.1 |
| 6 Operatives | 15.1 | 14.9 | 15.1 | 15.6 | 14.9 |
| 7 Others | 12.6 | 20.6 | 12.2 | 15.7 | 8.5 |
| ACTIVITY SECTOR | | | | | |
| 0 No | 2.4 | 18.9 | 0.3 | 0.4 | 0.4 |
| 1 Agricultural, Mining, Construction | 13.1 | 14.2 | 12.9 | 10.0 | 14.8 |
| 2 Manufacturing | 22.8 | 17.4 | 20.6 | 25.4 | 25.3 |
| 3 Transportation, Communication | 9.1 | 4.6 | 8.1 | 7.7 | 12.6 |
| 4 Wholesale and retail Trade | 14.9 | 17.4 | 14.9 | 16.5 | 13.1 |
| 5 Finance, Insurance | 4.3 | 2.5 | 4.0 | 5.1 | 4.6 |
| 6 Business and Personal Services | 9.1 | 12.8 | 9.4 | 7.9 | 8.4 |
| 7 Professional Services | 16.6 | 9.6 | 20.6 | 20.5 | 12.3 |
| 8 Public administration | 7.7 | 2.5 | 9.1 | 6.5 | 8.6 |

**Table 5: Frequencies of some qualitative variables (year 1992)**

Some simple results are visible in this table:
  - the proportion of white people is very weak in class A and significantly above the average in class D
  - a similar result is obtained with the highest level of education concerning these two classes
  - a special mention have to be made concerning the distribution of class C among the levels of education: the importance of the secondary level (high school without the final grade) and a post-graduate level slightly greater than the average are to be considered together with the great stability of this sub-population in a class constituted by older people with an important seniority in the same firm and moderate wages. This could be interpreted as a working population with a traditional qualification or on-the-job training implying a great stability but not a real career in their job.
- not surprisingly, managers and professional workers are merely present in class D, while class A is constituted mainly with unemployed and unqualified workers (Others are principally laborers)
- class C is made of subordinate occupations.

The classification used to represent the sectors of activity seems to be too aggregate to be conclusive[8].

## 6. THE TRAJECTORIES BETWEEN THE 4 MAIN CLASSES

Over the 2 507 individuals only 1 028 different trajectories are found, to be compared to the $4^9$ possible trajectories, it is clear that a trajectory cannot be conceived as a random process between the four classes.

Besides, 971 persons stay in the same class during 8 or 9 years over the 9 years observed. More precisely, only 34 persons (1.4%) stay in the main class A during the 9 years, 423 (17%) stay 8 or 9 years in class B, 218 (9%) stay in class C and 306 (12%) are present 8 or 9 years in class D. The rest of the sample, 1 536 persons (61%) are less stable but the changes observed are not numerous, of course less numerous than what have been measured on the classification in 7 classes. This is well establish with the probabilities computed with the same idea of a dominant position (table 6): 5 years or more in a class over the 9 years.

| Dominant Position | Size | Proba class A | Proba class B | Proba class C | Proba class D |
|---|---|---|---|---|---|
| A | 179 | **0.75** | 0.13 | 0.06 | 0.07 |
| B | 951 | 0.07 | **0.82** | 0.01 | 0.10 |
| C | 498 | 0.05 | 0.04 | **0.86** | 0.06 |
| D | 732 | 0.04 | 0.11 | 0.03 | **0.82** |

**Table 6: Probabilities to be one year in a class being most of the time in a given class**

A very small group appears to be in no dominant position over the whole period (table 7)

| No dominant position | Size | Proba class A | Proba class B | Proba class C | Proba class D |
|---|---|---|---|---|---|
|  | 147 | 0.34 | 0.33 | 0.13 | 0.29 |

**Table 7: Probabilities when no class has a dominant position**

---

[8] But using a more detailed classification is limited by the size of the sample used.

A very stable situation is visible and specially for class C. It may be seen on table 7 that the significant changes are only between A and B, B and D and D and B. For the little group with noticeable changes it is clear that class C is mostly isolated from the rest of the labor market.

This study may be continued with the measure of the frequencies computed over the transitions from a class to another one (table 8)

|      | AB   | AC   | AD   | BA   | BC   | BD   | CA   | CB   | CD   | DA   | DB   | DC   |
|------|------|------|------|------|------|------|------|------|------|------|------|------|
| Size | 554  | 177  | 242  | 492  | 159  | 1036 | 175  | 150  | 262  | 241  | 871  | 306  |
| %    | 0.12 | 0.04 | 0.05 | 0.11 | 0.03 | 0.22 | 0.04 | 0.03 | 0.06 | 0.05 | 0.19 | 0.07 |

**Table 8: Frequencies of transitions**

Again it is clear that the transitions occurred mainly between classes A, B et D. The number of improvements of the working situation (AB, BD) is close to the number of deterioration (BA, DB) and the moves from and to class C are very few.

It means that there is a clear difference between the hierarchy of the diverse situations observed and the dynamics of transitions: there is no reason for these transitions to follow the hierarchy, in a positive direction or the inverse. Class C is not a step towards the best state, class D.

The total number of the transitions showed in table 9 is 4665, that is a little less than 25% over the total number of possible transitions (8 x 2507), the major part is then constituted by the unchanged situation from one period to the next one.

These empirical probabilities may be used to build a transition matrix in order to model the changes of classes by a Markov chain. See in Table 9 the Markov transition matrix. For example, 0.57 is the probability to reach class A in the next year, starting from class A in the present year.

|   | A    | B    | C    | D    |
|---|------|------|------|------|
| A | 0.57 | 0.24 | 0.08 | 0.11 |
| B | 0.06 | 0.78 | 0.02 | 0.14 |
| C | 0.04 | 0.14 | 0.85 | 0.06 |
| D | 0.04 | 0.04 | 0.05 | 077  |

**Table 9: The Markov transition matrix**

This model relies on some important hypotheses concerning the factors influencing the transitions, precisely that the only factors which are taken into account for the constitution of the main classes are stable over the considered period,

Then it is possible to compute the stationary distribution[9] and to compare it to the distribution for the observed years. See Table 10.

|  | class A | class B | class C | class D |
|---|---|---|---|---|
| 1984 | .138 | .400 | .181 | .281 |
| 1988 | .110 | .381 | .199 | .309 |
| 1992 | .112 | .356 | .203 | .329 |
| stationary | .106 | .363 | .209 | .322 |

**Table 10: Distributions of three years and stationary distribution**

We can conclude that the situation observed, presented in table 10 for only 3 years of the period, is quite close from the limit situation computed with the Markov chain. This confirms the great stability of the main structure, even if a visible proportion of changes are observed mostly at the beginning and in the middle of the chain. The same conclusion is obtained when using the transition matrix produced with the classification in 7 classes.

We will continue the identification of the main factors influencing the transitions and determining the duration observed before these transitions.

## 7. CONCLUSION

A first of results establishes clearly the pertinence of a segmentation of the labor market in terms of intensity and stability of the occupation. A first half of the whole studied population exerts a very stable job, 60 % of this segment having a set of valuable conditions associated to this job (level and growth of the real wage).
Considering the trajectories followed by the individuals among the main possible situations some results are obtained:
    - a great stability is observed during the whole period studied, this means that most of the time each individual stays in the same situation from one year to the next one
    - the inferior part of what may be named the primary segment (class C) is the most stable situation
    - most of the observed transitions are between the two lower classes but a significant number of transitions is observed between these classes and the upper part of the primary segment.
A study of the factors influencing these different types of transitions would be valuable to understand the process of the trajectories observed.

---

[9] For an irreducible Markov chain (i.e. starting from any state it is possible to go to any state), the stationary distribution is defined as the invariant distribution, it is also the limit distribution when the time goes to infinity.

# APPENDIX 1

| | |
|---|---|
| AGEH | age of the head of the household |
| ANCH | number of years worked since the age of 18 |
| GRSALH | annual growth of hourly wage |
| HEXJH | annual work hours in extra jobs |
| HMJH | annual work hours in main job |
| HWMJH | annual number of weeks of work in main job |
| NBXJH | number of extra jobs |
| RSALH | real hourly wage |
| SENH | seniority in main jobs (in month) |
| SIZFAM | family size |
| VHWMJH | variation of work hours per week in main job, since 1982 |
| VWMJH | variation of weeks of work in main job, since 1982 |
| WMJH | number of weeks of work in main job |
| WOUTH | number of weeks out of labor force |
| WUNEH | number of weeks of unemployment |

**Table : Quantitative variables used in the Kohonen Classification**

# APPENDIX 2 (Principal Component Analysis)

**Graph 1**

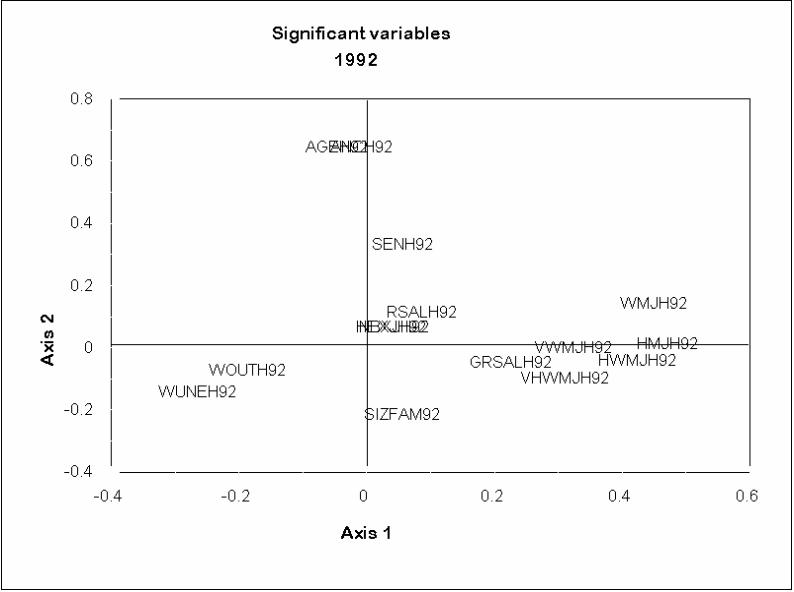

**Graph 2**

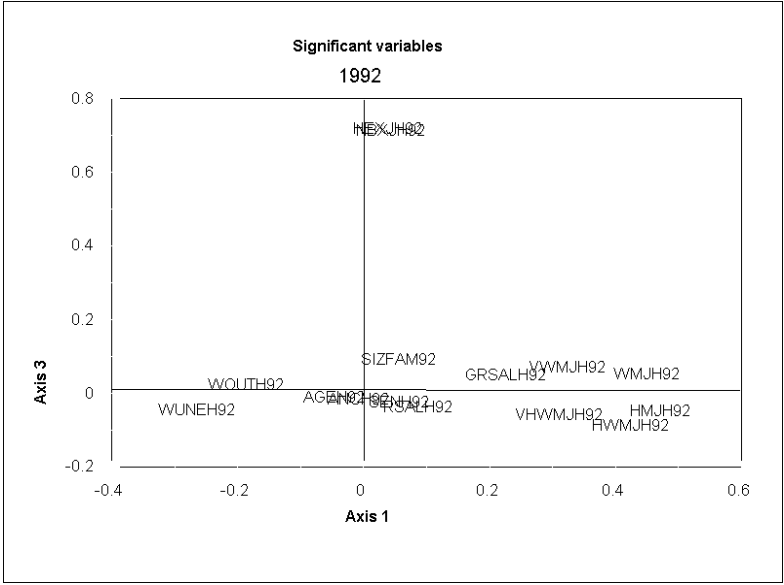

**Graph 3**

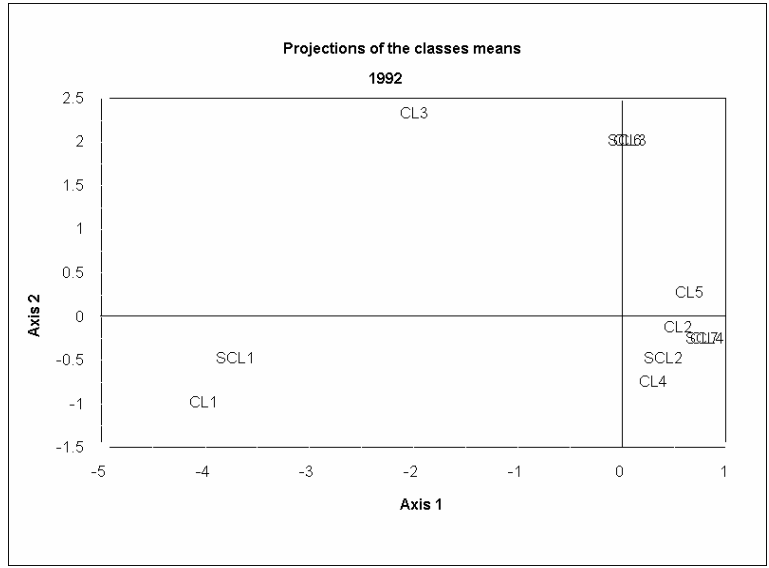

# APPENDIX 3

**Repartition of the qualitative variables for 1988 in the 4 main classes**

| in 1988 | Whole sample | Class A | Class B | Class C | Class D |
|---|---|---|---|---|---|
| RACE | | | | | |
| 1 Whites | 69.1 | 45.4 | 70.4 | 71.6 | 74.1 |
| 2 Blacks | 29.7 | 53.8 | 28.6 | 26.8 | 24.3 |
| EDUCATION | | | | | |
| 1 Primary | 0.9 | 1.8 | 0.2 | 2.0 | 0.6 |
| 2 Secondary | 18.9 | 33.1 | 12.2 | 29.0 | 15.7 |
| 3 Sec. achieved | 33.0 | 33.5 | 37.0 | 32.2 | 28.4 |
| 4 Post-sec. | 39.3 | 27.6 | 43.8 | 29.6 | 43.9 |
| 5 BA & more | 7.9 | 4.0 | 6.7 | 7.2 | 11.3 |
| OCCUPATION | | | | | |
| 0 No | 2.4 | 21.5 | 0.1 | 0 | 0 |
| 1-2 Managers, professionals | 35.1 | 14.6 | 34.6 | 33.0 | 44.2 |
| 4 Clerks | 11.7 | 12.0 | 13.0 | 14.0 | 8.6 |
| 5 Craftsmen | 17.9 | 8.0 | 19.2 | 20.0 | 18.4 |
| 6 Operatives | 14.1 | 13.8 | 14.1 | 14.8 | 13.8 |
| 7 Others | 13.5 | 26.9 | 13.8 | 13.8 | 8.2 |
| ACTIVITY SECTOR | | | | | |
| 0 No | 2.7 | 21.8 | 0.4 | 0.6 | 0.1 |
| 1 Agricultural, Mining, Construction | 13.7 | 14.5 | 14.0 | 10.0 | 15.3 |
| 2 Manufacturing | 23.6 | 16.4 | 22.0 | 28.8 | 24.7 |
| 3 Transportation, Communication | 8.8 | 4.0 | 7.6 | 8.6 | 12.0 |
| 4 Wholesale and retail Trade | 15.3 | 13.8 | 16.5 | 15.2 | 14.4 |
| 5 Finance, Insurance | 4.4 | 3.6 | 4.3 | 4.2 | 4.9 |
| 6 Business and Personal Services | 8.6 | 10.9 | 9.4 | 7.8 | 7.3 |
| 7 Professional Services | 15.4 | 12.4 | 17.3 | 17.6 | 12.8 |
| 8 Public administration | 7.5 | 2.5 | 8.5 | 7.2 | 8.4 |

**Frequencies of some qualitative variables (year 1988)**

**Repartition of the qualitative variables for 1984 in the 4 main classes**

| in 1984 | Whole sample | Class A | Class B | Class C | Class D |
|---|---|---|---|---|---|
| RACE | | | | | |
| 1 Whites | 67.2 | 49.1 | 69.6 | 70.9 | 70.2 |
| 2 Blacks | 29.9 | 48.6 | 27.6 | 26.0 | 26.4 |
| EDUCATION | | | | | |
| 1 Primary | 1.4 | 1.7 | 0.3 | 3.3 | 1.4 |
| 2 Secondary | 19.4 | 31.8 | 13.7 | 25.8 | 17.2 |
| 3 Sec. achieved | 39.0 | 43.9 | 41.6 | 35.8 | 34.9 |
| 4 Post-sec. | 34.7 | 20.5 | 40.1 | 28.0 | 38.4 |
| 5 BA & more | 5.5 | 2.0 | 4.2 | 7.1 | 8.1 |
| OCCUPATION | | | | | |
| 0 No | 4.1 | 28.6 | 0.1 | 0.2 | 0.3 |
| 1-2 Managers, professionals | 32.6 | 8.1 | 34.1 | 35.3 | 40.8 |
| 4 Clerks | 12.9 | 12.1 | 14.2 | 13.9 | 10.8 |
| 5 Craftsmen | 17.0 | 12.7 | 18.2 | 16.8 | 17.5 |
| 6 Operatives | 15.7 | 17.6 | 13.8 | 16.1 | 17.0 |
| 7 Others | 12.4 | 19.1 | 13.0 | 12.4 | 8.1 |
| ACTIVITY SECTOR | | | | | |
| 0 No | 4.7 | 29.2 | 0.9 | 0.9 | 0.6 |
| 1 Agricultural, Mining, Construction | 12.7 | 11.6 | 12.9 | 9.7 | 14.8 |
| 2 Manufacturing | 24.1 | 19.1 | 21.8 | 28.0 | 27.1 |
| 3 Transportation, Communication | 9.7 | 5.8 | 9.7 | 8.8 | 12.1 |
| 4 Wholesale and retail Trade | 15.4 | 11.6 | 17.1 | 15.5 | 14.9 |
| 5 Finance, Insurance | 4.1 | 1.7 | 4.3 | 5.1 | 4.5 |
| 6 Business and Personal Services | 7.6 | 10.4 | 7.4 | 7.1 | 7.0 |
| 7 Professional Services | 14.1 | 7.5 | 16.8 | 17.7 | 11.2 |
| 8 Public administration | 7.6 | 3.2 | 9.1 | 7.3 | 7.8 |

**Frequencies of some qualitative variables (year 1984)**